\documentclass[letterpaper, 10 pt, conference]{ieeeconf}
\IEEEoverridecommandlockouts
\usepackage{amsmath,amssymb,amsthm,mathtools}
\usepackage{bm}
\usepackage{mathrsfs}
\let\labelindent\relax
\usepackage{enumitem}
\usepackage{graphicx}
\usepackage[caption=false,font=footnotesize]{subfig}
\usepackage{xcolor}
\usepackage{cite}
\usepackage[hidelinks]{hyperref}
\usepackage[nameinlink]{cleveref}

\theoremstyle{plain}
\newtheorem{theorem}{Theorem}
\newtheorem{proposition}[theorem]{Proposition}
\newtheorem{lemma}[theorem]{Lemma}
\theoremstyle{definition}

\newtheorem{assumption}[theorem]{Assumption}
\newtheorem{remark}[theorem]{Remark}

\crefname{equation}{}{}
\Crefname{equation}{Eq.}{Eqs.}
\crefname{theorem}{Theorem}{Theorems}
\Crefname{theorem}{Theorem}{Theorems}
\crefname{lemma}{Lemma}{Lemmas}
\Crefname{lemma}{Lemma}{Lemmas}
\crefname{proposition}{Proposition}{Propositions}
\Crefname{proposition}{Proposition}{Propositions}
\crefname{definition}{Definition}{Definitions}
\Crefname{definition}{Definition}{Definitions}
\crefname{assumption}{Assumption}{Assumptions}
\Crefname{assumption}{Assumption}{Assumptions}
\crefname{remark}{Remark}{Remarks}
\Crefname{remark}{Remark}{Remarks}
\crefname{figure}{Fig.}{Figs.}
\Crefname{figure}{Fig.}{Figs.}
\crefname{table}{Table}{Tables}
\Crefname{table}{Table}{Tables}
\crefname{section}{Section}{Sections}
\Crefname{section}{Section}{Sections}

\newcommand{\Def}{\coloneqq}
\newcommand{\RR}{\mathbb{R}}
\newcommand{\NN}{\mathbb{N}}
\newcommand{\eps}{\varepsilon}
\newcommand{\timehor}{\bm T_h}
\newcommand{\sa}{\mathsf a}
\newcommand{\VV}{\bm V}
\newcommand{\nv}{\mathsf n}
\newcommand{\tup}[1]{\textup{(}#1\textup{)}}

\newcommand{\Sh}{S_{h}}
\newcommand{\Tdt}{\mathcal T_{h}}
\newcommand{\QT}{\mathcal Q_{T}}
\newcommand{\Bnd}{\mathscr B_0}
\newcommand{\sub}{\mathfrak S}
\DeclarePairedDelimiter{\abs}{\lvert}{\rvert}
\DeclarePairedDelimiter{\norm}{\lVert}{\rVert}
\DeclarePairedDelimiter{\set}{\{}{\}}
\DeclarePairedDelimiter{\ceil}{\lceil}{\rceil}
\DeclareMathOperator{\dist}{dist}

\title{\LARGE \bf A Derivative-Free Bellman Equation for Dynamic
Time-to-Collision:\\ Two-Sided Certificates and Explicit ReLU
Approximation}

\author{Hossein Nick Zinat Matin, IEEE member%
\thanks{$^{1}$H. Nick Zinat Matin is with LIX, \'Ecole Polytechnique, and IP Paris,
France, {\tt\small hossein.matin@polytechnique.edu}.}%
\thanks{§This work benefited from Hi! PARIS and State funding managed by the French National Research Agency
(ANR) under the France 2030 program, reference ANR-23-IACL-0005.}
}

\begin{document}
\maketitle
\thispagestyle{empty}
\pagestyle{empty}

\begin{abstract}
Dynamic time-to-collision (TTC) is the first contact time of a
prediction flow and defines a feedback map implicitly. We study its
capped value function $\VV$ for Carath\'eodory dynamics with bounded,
measurable forcing. A self-contained argument establishes local
Lipschitz continuity on the transversal collision domain. A
derivative-free stopped-step Bellman equation then characterizes
$\VV$ uniquely: it advances the prediction for a fixed duration or
until contact and adds the elapsed time. A two-sided error bound uses
only the candidate's residual along complete stopped trajectories and
its error at their endpoints. Exact continuum iteration recovers
$\VV$ after finitely many sweeps, including one additional sweep for
initializations that violate the boundary condition. We also
characterize an optional-duration infimum operator: its fixed points
are the dynamic-programming subsolutions, so its residual alone does
not identify $\VV$. Finally, on compact subsets of the transversal
domain, an optimization-free ReLU construction has explicit size,
controlled local residual, and robustness to errors in its nodal TTC
values. The optimal velocity model provides the running application.
\end{abstract}

\section{Introduction}

Time-to-collision (TTC) estimates the first contact time under
assumptions on future motion \cite{hayward1971near,lee1976theory, matin2025second}.
The dynamic-based TTC of \cite{matin2026TTCandHJB} predicts contact
from a car-following flow and enters a feedback extension of the
Optimal Velocity Model (OVM) \cite{matin2020nonlinear, matin2020nonlinear2}. This makes the prediction consistent with the model used to define the feedback (cf. \cite{matin2024analytical, nick2022near}). Evaluating the map requires forward integration, motivating an explicit approximation
for repeated evaluation.

In \cite{matin2026TTCandHJB}, the TTC is identified with the value
function $\VV$ of a minimum-time Hamilton--Jacobi reachability problem
\cite{mitchell2005time,bansal2017hamilton} and approximated by a
physics-informed network \cite{raissi2019}. Its differential HJB
residual and a boundary loss bound the $L^2$ approximation error for
$C^1$ candidates; a reverse loss bound is not established. Bounded
forcing alone does not prevent a zero differential residual, but
additional jump and nondegeneracy conditions can obstruct an exact
$C^1$ solution. That work also gives an integrated characteristic
identity. Here we study a fixed stopped-step identity and separately
compare it with an optional-duration infimum operator, introduced
below as an alternative formulation.

This paper studies the feedback map $\VV$ through a different
derivative-free equation, the \emph{stopped-step} Bellman equation:
advance the prediction for a fixed duration $h$, or until contact,
whichever comes first, and add the elapsed time. The equation is
evaluated along the prediction flow, needs no derivative of the
candidate, and has $\VV$ as its unique solution with a two-sided
residual certificate. Our contributions are:
\begin{itemize}[leftmargin=*,nosep]
\item A self-contained regularity theory for the first-contact value
function of a Carath\'eodory flow with transversal contact
(\cref{prop:DPP,prop:regularity}, Appendix), allowing merely bounded
time-dependent inputs and models defined only up to contact.
\item The stopped-step equation: uniqueness of its solution, a
two-sided a priori error bound in terms of the residual along
trajectories, and exact recovery of $\VV$ after finitely many sweeps
from any initialization (\cref{thm:certificate,prop:iteration}). The
certificate permits arbitrary finite-valued candidates on the
trajectories and their reached endpoints.
\item A characterization of the optional-duration form: its
fixed points are exactly the dynamic-programming subsolutions, a
convex family with maximum $\VV$, so its residual alone cannot
identify the TTC (\cref{prop:envelope}).
\item An optimization-free ReLU realization of $\VV$ with explicit
size, controlled stopped-step residual, and robustness to errors in
the nodal values (\cref{thm:relu}).
\end{itemize}
The underlying mechanisms are classical: minimum-time regularity
with measurable time dependence \cite{veliov1997lipschitz},
maximal-subsolution characterizations \cite{motta2014asymptotic},
linear-programming and discrete dynamic programming
\cite{manne1960linear,defarias2003linear,bardi1990approximation}, and
ReLU interpolation \cite{yarotsky2017,he2020relu}. Oustry and Tacchi
\cite{oustry2023minimal} optimize differentiable HJB subsolutions for
controlled dynamics with Lipschitz data. Liu \cite{liu2026verifiable}
derives verifiable error bounds from differential residuals of
$C^1$ candidates for Lyapunov and HJB equations. Our narrower setting
is a prescribed flow with finite first-contact chains: it permits a
simple integrated residual certificate for nondifferentiable
candidates and bounded measurable forcing. The contribution lies in combining these mechanisms for dynamic TTC
under the stated low-regularity assumptions.

\section{Prediction Dynamics and Dynamic-Based TTC}
\label{sec:model}

\subsection{Closed-loop model}

The controlled system has state $\phi\in\RR^p$ and dynamics
\begin{equation}\label{E:closed-loop}
\dot\phi=\mathsf F\bigl(t,\phi,\tau(t,\phi)\bigr),\qquad
\tau(t,\phi)\Def\VV\bigl(t,\pi(\phi)\bigr),
\end{equation}
where $\tau$ is the dynamic-based TTC, $\VV$ is the value function
defined in \cref{E:value-function} below from a prediction model, and
$\pi:\RR^p\to\RR^d$ extracts from the actual state the initial
condition of the prediction. Thus $\tau$ is a time-varying static
state feedback, defined implicitly through the prediction model, as a
model predictive controller is; this paper concerns the regularity
and explicit approximation of this map. No closed-loop
well-posedness assumption is needed for our results.
For the car-following instance, \cite[Thm.~3.38 and Rem.~3.41]{matin2026TTCandHJB}
establishes existence and excludes contacts with negative relative
velocity under its assumptions, with additional conditions for
uniqueness and terminal extension.
A lower bound on $\VV$ concerns the first contact time
\emph{under the prediction model}; it does not establish closed-loop
collision avoidance.

\textbf{Car-following instance.} A platoon in which vehicle
$n\in\set{1,\dots,N}$, with position $x_n$ and velocity $v_n$, follows
vehicle $n-1$; vehicle $0$ is an exogenous leader with
$\dot v_{\mathsf L}=\sa_{\mathsf L}$, where
\begin{equation}\label{E:accel_L}
\sa_{\mathsf L}\in L^\infty(0,\timehor),\quad
0\le v_{{\mathsf L},\circ}+\!\int_0^t\!\sa_{\mathsf L}\,dr\le v_{\max},
\ \ t\in[0,\timehor],
\end{equation}
$v_{{\mathsf L},\circ}$ being the initial leader velocity and
$v_{\max}$ the maximal speed. The TTC-feedback OVM of
\cite{matin2026TTCandHJB} reads
\begin{equation}\label{E:main_dynamics}
\dot x_n=v_n,\qquad
\dot v_n=\alpha\bigl\{\mathcal V(x_{n-1}-x_n)-v_n\bigr\}
+\beta\,\frac{v_{n-1}-v_n}{\tau_n^{2}},
\end{equation}
with $\alpha,\beta>0$,
$\mathcal V(r)\Def\tanh(r-2)+\tanh 2$ the smooth, bounded, increasing,
$1$-Lipschitz optimal velocity function
\cite{nick2022near,bando1995dynamical}, and
$\tau_n(t)=\VV(t,\pi_n(\phi(t)))$ with
$\pi_n(\phi)\Def(x_{n-1}-x_n,\,v_{n-1}-v_n,\,v_{n-1})$. The
prediction for the pair $(n-1,n)$ is conditional on the acceleration
of vehicle $n-1$ over the horizon, prescribed as in \cref{E:accel_L};
for a single leader--follower pair this is the leader input
\cite{matin2024analytical,matin2026TTCandHJB}.
All values and Bellman successors below use the same prescribed
forecast. Updating that forecast changes the map and requires a
separate consistency analysis. The fixed terminal time also gives
$\VV(\timehor,\cdot)=0$; using \cref{E:main_dynamics} at that time
requires a separate terminal or gain-regularization convention.
The approximation results do not supply such a controller design.

\subsection{Prediction dynamics}

Let $U\subset\RR^d$ be open (the state space of the prediction model),
$g\in C^1(U)$ the \emph{collision margin},
$\mathcal Q\Def\set{\bm y\in U:g(\bm y)\le0}$ the collision set, closed
relative to $U$, and $P\Def\set{\bm y\in U:g(\bm y)\ge0}$ the region in
which the model is physically well defined. The prediction model is
\begin{equation}\label{E:general-field}
\dot{\bm y}=\mathsf G(t,\bm y),\qquad
\mathsf G:[0,\timehor]\times P\to\RR^d .
\end{equation}
We write $z=(s,\bm y)$ for a point of $[0,\timehor]\times U$,
$m\Def d+1$, and $\abs{\cdot}$ for Euclidean norms. A \emph{solution
in $P$} from $(s,\bm y)\in[0,\timehor]\times P$ is an absolutely
continuous $x:[s,s+\sigma]\to P$, $0\le\sigma\le\timehor-s$, with
$x(s)=\bm y$ and $\dot x=\mathsf G(t,x)$ a.e.; it is \emph{maximal} if
it is not the restriction of a solution in $P$ on a longer interval.

\begin{assumption}\label{ass:field}
\begin{enumerate}[label=\tup{A\arabic*},leftmargin=*,nosep]
\item\label{A1} $\mathsf G$ is measurable in $t$, and locally bounded
and locally Lipschitz in $\bm y$ on $P$, up to the collision surface,
with constants independent of $t$: for every compact $\bm C\subset P$
there are $M_{\bm C},L_{\bm C}$ such that, for a.e.\ $t$ and all
$\bm y,\bm y'\in\bm C$, $\abs{\mathsf G(t,\bm y)}\le M_{\bm C}$ and
$\abs{\mathsf G(t,\bm y)-\mathsf G(t,\bm y')}\le L_{\bm C}\abs{\bm y-\bm y'}$.
\item\label{A2} Forward completeness up to the boundary of $P$: from
every $(s,\bm y)\in[0,\timehor]\times P$ there is a maximal solution
in $P$, defined on a compact interval $[s,s+\sigma(s,\bm y)]$, so that
no solution escapes to infinity, or to $\partial U$, before reaching
$\timehor$ or leaving $P$ \tup{if $s+\sigma<\timehor$ then
$g(x(s+\sigma))=0$, since otherwise Carath\'eodory's existence theorem
on a closed ball in $\set{g>0}$ would continue the solution}. For $g(\bm y)>0$ the first contact time
$\varsigma(s,\bm y)\le\sigma(s,\bm y)$ is the first time the solution
reaches $\mathcal Q$, or $\timehor-s$ if there is none, and the
prediction is the solution on $[s,s+\varsigma(s,\bm y)]$.
\item\label{A3} The normal velocity
$\nv(t,\bm y)\Def\nabla g(\bm y)\cdot\mathsf G(t,\bm y)$ is continuous
on $[0,\timehor]\times P$.
\end{enumerate}
\end{assumption}

By \ref{A1} and Gronwall's inequality, two solutions in $P$ from the
same datum coincide on their common interval
\cite[Ch.~2, \S1]{filippov1988differential}. For $g(\bm y)>0$ we
denote by $\Phi(\kappa;s,\bm y)$, $0\le\kappa\le\varsigma(s,\bm y)$, the
state of the prediction at time $s+\kappa$. On this domain
$\kappa\mapsto\Phi(\kappa;s,\bm y)$ is Lipschitz, and the semigroup property
$\Phi(\kappa_2;s+\kappa_1,\Phi(\kappa_1;s,\bm y))
=\Phi(\kappa_1+\kappa_2;s,\bm y)$ holds whenever
$\kappa_1+\kappa_2\le\varsigma(s,\bm y)$. Nothing beyond contact is defined
or used: the operator of \cref{sec:stopped} evaluates the flow only on
$[0,\varsigma]$, and the proof of \cref{prop:regularity} uses a local
extension of $\mathsf G$ that does not affect contact times. Condition
\ref{A3} says that the rate at which $g$ changes along trajectories
is continuous even though $\mathsf G$ is only measurable in $t$; it
holds for second-order car-following models with Lipschitz
acceleration laws and bounded inputs, where $g$ is the gap and $\nv$
the relative velocity. A uniform version of \ref{A2}, trajectories
from a compact set of initial data staying in a compact subset of $P$,
follows from \ref{A1}--\ref{A2} (\cref{lem:uniform-bound} in the
Appendix).

\textbf{The OVM instance.} For the car-following application,
\cref{E:general-field} is the OVM, i.e., \cref{E:main_dynamics}
without the $\beta$-term, started from $\pi_n$ of the current state.
Denoting the predicted positions and velocities of leader and
follower by $(\xi_{\mathsf L},\zeta_{\mathsf L})$ and $(\xi,\zeta)$, the
predicted gap $\Xi\Def\xi_{\mathsf L}-\xi$ and relative velocity
$\Lambda\Def\zeta_{\mathsf L}-\zeta$, the state is
$\bm y=(\Xi,\Lambda,\zeta_{\mathsf L})\in\RR^3$, $U=\RR^3$,
$g(\bm y)=\Xi$, $P=\set{\Xi\ge0}$, and
\begin{equation}\label{E:new-accel}
\begin{aligned}
\dot\Xi&=\Lambda,\qquad \dot\zeta_{\mathsf L}=\sa_{\mathsf L}(t),\\
\dot\Lambda&=\sa_{\mathsf L}(t)-\alpha\bigl\{\mathcal V(\Xi)
-\zeta_{\mathsf L}+\Lambda\bigr\}.
\end{aligned}
\end{equation}
\Cref{ass:field} holds: the forcing $\sa_{\mathsf L}$ enters
additively and cancels in $\mathsf G(t,\bm y)-\mathsf G(t,\bm y')$, so
the field is globally Lipschitz in $\bm y$ with a constant
$1+2\alpha$ independent of $t$ and, since $\abs{\mathsf G(t,0)}$ is
bounded, of linear growth, which gives \ref{A2} by Gronwall's
inequality; and $\nv=\Lambda$ is a state component, hence continuous.
Here $m=4$.

\subsection{Dynamic-based TTC and the collision domain}

The TTC is the value function of the first-contact problem
\begin{equation}\label{E:value-function}
\VV(s,\bm y)\Def\min\Bigl\{\timehor-s,\;
\inf_{\substack{0\le\kappa\le\varsigma(s,\bm y)\\
\Phi(\kappa;s,\bm y)\in\mathcal Q}}\kappa\Bigr\},
\end{equation}
with $\inf\varnothing=+\infty$, i.e.\ $\VV(s,\bm y)=\varsigma(s,\bm y)$ for
$g(\bm y)>0$, and $\VV\Def0$ on $\QT\Def[0,\timehor]\times\mathcal Q$.
Thus $0\le\VV\le\timehor-s$, with equality on the right exactly when
no contact is predicted before $\timehor$. A tangential touch, $g=0$
with $\nv=0$, counts as a contact: $\VV$ is the time until the
predicted margin first closes, whatever happens afterwards. Contacts are sought in the
\emph{transversal collision domain}
\begin{equation}\label{E:mathbb-D}
\begin{split}
\mathbb D\Def\bigl\{(s,\bm y)\in[0,\timehor)\times U:{}&
g(\bm y)>0,\ \VV(s,\bm y)<\timehor-s,\\
&\nv\bigl(s+\VV,\Phi(\VV;s,\bm y)\bigr)<0\bigr\},
\end{split}
\end{equation}
the states whose first contact with $\set{g=0}$ is transversal. Writing
$\Lambda(\kappa;s,\bm y)$ for the $\Lambda$-component of
$\Phi(\kappa;s,\bm y)$, for the OVM the last condition reads
$\Lambda(\VV;s,\bm y)<0$: the predicted relative velocity is negative
at impact. For $(s,\bm y)\in\mathbb D$ we call
$(s+\VV,\Phi(\VV;s,\bm y))\in\QT$ the \emph{impact point}. We also use
the stopped flow $\Phi^{\mathcal Q}(\kappa;s,\bm y)\Def\Phi(\min\{\kappa,
\VV(s,\bm y)\};s,\bm y)$ for $g(\bm y)>0$ and
$\Phi^{\mathcal Q}(\kappa;s,\bm y)\Def\bm y$ on $\mathcal Q$, which
remains at the impact point once $\mathcal Q$ is reached, and
$\Omega\Def\mathbb D\cup\QT$.

\begin{proposition}[Dynamic programming]\label{prop:DPP}
Let $(s,\bm y)\in\mathbb D$ and $\nu\Def\VV(s,\bm y)$. For every
$\kappa\in[0,\nu)$, $(s+\kappa,\Phi(\kappa;s,\bm y))\in\mathbb D$ and
$\VV(s+\kappa,\Phi(\kappa;s,\bm y))=\nu-\kappa$.
\end{proposition}

\begin{proof}
Let $z'\Def(s+\kappa,\Phi(\kappa;s,\bm y))$. By the semigroup
property, $\Phi(u;z')=\Phi(u+\kappa;s,\bm y)$ for $0\le u\le\nu-\kappa$.
Since $\nu$ is the first contact time from $(s,\bm y)$ and
$\kappa<\nu$, we have $g(\Phi(\kappa;s,\bm y))>0$, no contact from $z'$
before $u=\nu-\kappa$, and a contact at
$u=\nu-\kappa<\timehor-(s+\kappa)$. Hence $\VV(z')=\nu-\kappa$, and the
impact point of $z'$ is that of $(s,\bm y)$, where $\nv<0$. Therefore
$z'\in\mathbb D$.
\end{proof}

\begin{proposition}[Regularity]\label{prop:regularity}
$\mathbb D$ is relatively open in $[0,\timehor)\times U$, and $\VV$ is
locally Lipschitz on $\mathbb D$, hence Lipschitz on every compact
$\bm K\subset\mathbb D$, with a constant denoted $L_{\VV}(\bm K)$.
Moreover $(\kappa,z)\mapsto\Phi(\kappa;z)$ is continuous on
$\set{(\kappa,z):z\in\mathbb D,\ 0\le\kappa\le\VV(z)}$.
\end{proposition}

\begin{proof}
Fix $z_0=(s_0,\bm y_0)\in\mathbb D$, $\nu_0\Def\VV(z_0)>0$, the impact
state $q\Def\Phi(\nu_0;z_0)$, the impact time $t_*\Def s_0+\nu_0$, and
$a\Def-\nv(t_*,q)>0$.

\emph{Step 0: local extension.} The pre-contact path
$\Gamma\Def\Phi([0,\nu_0];z_0)$ is a compact subset of $U$. Choose
$\eta>0$ such that the closed $\eta$-neighborhood $\Gamma_\eta$ of
$\Gamma$ is a compact subset of $U$, and set $\bm C\Def\Gamma_\eta\cap P$, a
compact subset of $P$. Extend each component of $\mathsf G(t,\cdot)$
from $\bm C$ to $\RR^d$ by the McShane formula with the constant
$L_{\bm C}$ of \ref{A1}; on the null set of $t$ where the inequalities
of \ref{A1} fail, set the extension to $0$. The extended field
$\tilde{\mathsf G}$ agrees with $\mathsf G$ on $\bm C$ for a.e.\ $t$,
is measurable in $t$ \tup{infimum over a countable dense subset of
$\bm C$}, globally Lipschitz in $\bm y$ with constant
$\tilde L\le\sqrt d\,L_{\bm C}$, and bounded by $M_{\bm C}+\tilde L
\abs{\bm y-q}$ near $q$; it has linear growth, so its flow
$\tilde\Phi$ is defined on $[s,\timehor]$ for every datum and is
jointly continuous. By uniqueness, $\tilde\Phi(\kappa;z)=\Phi(\kappa;z)$
whenever the latter stays in $\bm C$; in particular on $[0,\nu_0]$ for
$z=z_0$. Set $F(z,\kappa)\Def g(\tilde\Phi(\kappa;z))$ wherever
$\tilde\Phi(\kappa;z)\in U$.

\emph{Step 1: Lipschitz dependence.} Choose $\rho,\eps_0>0$ with
$s_0+\rho<\timehor$, $2\rho\le\eps_0$ and
$\nu_0+2\eps_0\le\timehor-s_0-\rho$, let
$\mathcal N\Def([s_0-\rho,s_0+\rho]\cap[0,\timehor))\times\bar B(\bm y_0,\rho)$
and $\bar\kappa\Def\nu_0+\eps_0$, so that
$\bar\kappa+\eps_0\le\timehor-s$ on $\mathcal N$. By joint continuity of
$\tilde\Phi$, after shrinking $\rho$ and $\eps_0$ \tup{keeping
$2\rho\le\eps_0$}, the points
$\tilde\Phi(\kappa;z)$, $z\in \mathcal N$, $0\le\kappa\le\bar\kappa+\eps_0$, lie
in the interior of $\Gamma_\eta$, where $g$ is Lipschitz with a constant
$L_g$, and $\abs{\tilde{\mathsf G}}\le\tilde M$ there. Then
$\abs{\tilde\Phi(\kappa;z)-\tilde\Phi(\kappa';z)}\le\tilde M
\abs{\kappa-\kappa'}$, Gronwall's inequality gives
$\abs{\tilde\Phi(\kappa;s,\bm y)-\tilde\Phi(\kappa;s,\bm y')}\le
e^{\tilde L\bar\kappa}\abs{\bm y-\bm y'}$ for $\kappa\le\bar\kappa$,
and, for $s'=s+\delta$ with both $(s,\bm y),(s',\bm y)\in \mathcal N$, the bound
$\bar\kappa+\eps_0\le\timehor-s'$ and $\delta\le2\rho\le\eps_0$ give
$\kappa+\delta<\timehor-s$ and $\kappa+\delta\le\bar\kappa+\eps_0$, so
the semigroup identity
$\tilde\Phi(\kappa+\delta;s,\bm y)=\tilde\Phi(\kappa;s',
\tilde\Phi(\delta;s,\bm y))$ and Gronwall on $[0,\kappa]$ yield
$\abs{\tilde\Phi(\kappa;s',\bm y)-\tilde\Phi(\kappa;s,\bm y)}
\le(e^{\tilde L\bar\kappa}+1)\tilde M\delta$. Since $\mathcal N$ is a product,
$F$ is Lipschitz on $\mathcal N\times[0,\bar\kappa]$ with a constant $L_F$.

\emph{Step 2: strict decrease near impact.} For $x\in U$ and
a.e.\ $t$,
\begin{equation*}
\begin{split}
\nabla g(x)\cdot\tilde{\mathsf G}(t,x)-\nv(t,q)
={}&[\nabla g(x)-\nabla g(q)]\cdot\tilde{\mathsf G}(t,x)\\
&+\nabla g(q)\cdot[\tilde{\mathsf G}(t,x)-\mathsf G(t,q)],
\end{split}
\end{equation*}
using $\tilde{\mathsf G}(t,q)=\mathsf G(t,q)$ for a.e.\ $t$. By
continuity of $\nabla g$ at $q$, the global Lipschitz bound of
$\tilde{\mathsf G}$ and its local boundedness, there is $\rho_1>0$
with $\bar B(q,\rho_1)\subset U$ such that the right-hand side is at
most $a/4$ in absolute value for a.e.\ $t$ and all
$\abs{x-q}\le\rho_1$. By \ref{A3} there is $\eps_1\in(0,\timehor-t_*]$
with $\nv(t,q)\le-3a/4$ for $\abs{t-t_*}\le\eps_1$; shrinking
$\rho_1,\eps_1$, also $\nv<0$ on
$[t_*-\eps_1,t_*+\eps_1]\times(\bar B(q,\rho_1)\cap P)$ by continuity
of $\nv$. Hence
\begin{equation}\label{eq:ae-decrease}
\begin{gathered}
\nabla g(x)\cdot\tilde{\mathsf G}(t,x)\le-c\Def-a/2,\\
\text{for a.e.\ }t\in[t_*-\eps_1,t_*+\eps_1],\quad\abs{x-q}\le\rho_1.
\end{gathered}
\end{equation}
Fix $\eps\in(0,\min\{\eps_0,\eps_1/2,\nu_0,\rho_1/(2\tilde M)\}]$, so
that $\tilde\Phi(\kappa;z_0)\in\bar B(q,\rho_1/2)$ for
$\abs{\kappa-\nu_0}\le\eps$ by Step~1, and shrink $\mathcal N$ to a
compact relative neighborhood $\mathcal N'$ of $z_0$, of $s$-radius at
most $\eps_1/2$, such that $\tilde\Phi(\kappa;z)\in\bar B(q,\rho_1)$
for $z\in \mathcal N'$ and $\kappa\in[\nu_0-\eps,\nu_0+\eps]$ \tup{joint
continuity}. For such $z=(s,\bm y)$ and $\kappa$,
$\abs{s+\kappa-t_*}\le\eps_1$, so \cref{eq:ae-decrease} applies along
the trajectory. Since
$\kappa\mapsto g(\tilde\Phi(\kappa;z))$ is absolutely continuous with
derivative $\nabla g(\tilde\Phi)\cdot\tilde{\mathsf G}(s+\kappa,
\tilde\Phi)$ a.e., \cref{eq:ae-decrease} integrates to
\begin{equation}\label{eq:decrease}
\begin{gathered}
F(z,\kappa_2)-F(z,\kappa_1)\le-c(\kappa_2-\kappa_1),\\
\nu_0-\eps\le\kappa_1\le\kappa_2\le\nu_0+\eps,
\end{gathered}
\end{equation}
for all $z\in \mathcal N'$.

\emph{Step 3: the first contact time.} Since $F(z_0,\cdot)>0$ on the
compact interval $[0,\nu_0-\eps]$ and $F$ is Lipschitz, after
shrinking $\mathcal N'$ we have $F(z,\kappa)>0$ for $z\in \mathcal N'$,
$\kappa\in[0,\nu_0-\eps]$. By \cref{eq:decrease},
$F(z_0,\nu_0+\eps)\le-c\eps$, so after shrinking $\mathcal N'$ once more,
$F(z,\nu_0+\eps)<0$ for $z\in \mathcal N'$. Hence $F(z,\cdot)$ has exactly one
zero $\kappa^*(z)\in(\nu_0-\eps,\nu_0+\eps)$, is positive before it,
and $\kappa^*(z)<\bar\kappa<\timehor-s$. Before $\kappa^*(z)$ the
extended trajectory lies in $P\cap \Gamma_\eta=\bm C$, where
$\tilde{\mathsf G}=\mathsf G$, so by uniqueness it is the prediction
from $z$, and $\VV(z)=\varsigma(z)=\kappa^*(z)$; in particular
$\Phi(\kappa;z)=\tilde\Phi(\kappa;z)$ for $z\in\mathcal N'$,
$0\le\kappa\le\VV(z)$, which gives the continuity claim. Moreover $g(\bm y)=F(z,0)>0$,
and the impact point of $z$ lies in $\bar B(q,\rho_1)$ at a time in
$[t_*-\eps_1,t_*+\eps_1]$, where $\nv<0$ by the choice of
$\rho_1,\eps_1$. Hence $\mathcal N'\subset\mathbb D$,
which proves relative openness.

\emph{Step 4: Lipschitz bound.} For $z,z'\in \mathcal N'$, \cref{eq:decrease}
gives $\abs{F(z',\kappa^*(z))}=\abs{F(z',\kappa^*(z))-F(z',\kappa^*(z'))}
\ge c\abs{\kappa^*(z)-\kappa^*(z')}$, while Step 1 gives
$\abs{F(z',\kappa^*(z))}=\abs{F(z',\kappa^*(z))-F(z,\kappa^*(z))}\le
L_F\abs{z-z'}$. Hence $\abs{\VV(z)-\VV(z')}\le(L_F/c)\abs{z-z'}$.
A locally Lipschitz function is Lipschitz on compact subsets, by a
standard covering argument.
\end{proof}

In particular $\VV$ is continuous on $\mathbb D$. Since $\mathbb D$
is relatively open and $\QT$ relatively closed in
$[0,\timehor]\times U$, both are Borel and $\VV$ is Borel on $\Omega$.
The Lipschitz constants need not remain uniform near grazing,
where $\VV$ can be discontinuous; the resulting closed-loop issues
are analyzed in \cite{matin2026TTCandHJB}.

\section{The Stopped-Step Bellman Equation}
\label{sec:stopped}

Let $\Bnd$ be the set of bounded functions $W:\Omega\to\RR$ with
$W=0$ on $\QT$; functions are evaluated pointwise along trajectories,
so we work with bounded functions rather than $L^\infty$ classes, and
$\VV\in\Bnd$ pointwise. Fix a step $h>0$. For $z=(s,\bm y)\in\mathbb D$
let
\begin{equation}\label{E:ell}
\ell_h(z)\Def\min\{h,\VV(z)\},\qquad
z_\kappa\Def\bigl(s+\kappa,\Phi^{\mathcal Q}(\kappa;z)\bigr),
\end{equation}
and define, for a bounded $W:\Omega\to\RR$, $\Sh[W]\Def0$ on $\QT$ and
\begin{equation}\label{E:S-operator}
\Sh[W](z)\Def\ell_h(z)+W\bigl(z_{\ell_h(z)}\bigr),\qquad z\in\mathbb D .
\end{equation}
Operationally, \cref{E:S-operator} integrates the prediction for
duration $h$, or until contact if that occurs first, and adds the
elapsed time; no knowledge of $\VV$ is needed beyond event detection,
and no derivative of $W$ enters. The evaluation point $z_{\ell_h(z)}$
lies in $\Omega$: in $\mathbb D$ if $h<\VV(z)$ by \cref{prop:DPP}, and
at the impact point, in $\QT$, otherwise. Hence $\Sh$ maps bounded
functions on $\Omega$ to bounded functions and $\Bnd$ into itself; it
is monotone, $W_1\le W_2\Rightarrow\Sh[W_1]\le\Sh[W_2]$, and
non-expansive,
\begin{equation}\label{eq:non-expansion}
\abs[\big]{\Sh[W_1]-\Sh[W_2]}(z)\le\abs{W_1-W_2}\bigl(z_{\ell_h(z)}\bigr),
\quad z\in\mathbb D .
\end{equation}
\begin{lemma}[Fixed point]\label{lem:fixed-point}
$\Sh[\VV]=\VV$ on $\Omega$ for every $h>0$.
\end{lemma}

\begin{proof}
On $\QT$ both sides vanish. For $z\in\mathbb D$ with
$\nu\Def\VV(z)$: if $h<\nu$, \cref{prop:DPP} gives
$\VV(z_h)=\nu-h$, so $\Sh[\VV](z)=h+\nu-h=\nu$; if $h\ge\nu$, then
$\ell_h(z)=\nu$, $z_\nu$ is the impact point, $\VV(z_\nu)=0$, and
$\Sh[\VV](z)=\nu$.
\end{proof}

\begin{theorem}[Uniqueness and two-sided certificate]
\label{thm:certificate}
Let $\bm A\subset\mathbb D$ contain the trajectories of its points up
to impact, i.e.\ $(s+\kappa,\Phi(\kappa;s,\bm y))\in\bm A$ whenever
$(s,\bm y)\in\bm A$ and $0\le\kappa<\VV(s,\bm y)$, and let
$\Gamma_{\bm A}\Def\set{z_{\VV(z)}:z\in\bm A}$ be its reached impact
set. Let $W:\bm A\cup\Gamma_{\bm A}\to\RR$ be finite valued with
$\abs{W}\le b$ on $\Gamma_{\bm A}$ and
\begin{equation}\label{eq:residual}
\abs[\big]{W-\Sh[W]}\le\delta\quad\text{on }\bm A .
\end{equation}
Here \cref{E:S-operator} defines $\Sh[W]$ on $\bm A$ using only
values in $\bm A\cup\Gamma_{\bm A}$.
Then, for all $z\in\bm A$,
\begin{equation}\label{eq:two-sided}
\abs{W(z)-\VV(z)}\le\ceil[\big]{\VV(z)/h}\,\delta+b .
\end{equation}
In particular, $\VV$ is the unique fixed point of $\Sh$ in $\Bnd$.
\end{theorem}

\begin{proof}
Fix $z\in\bm A$, $\nu\Def\VV(z)$, $N\Def\ceil{\nu/h}$, and define
$z_0\Def z$, $z_{j+1}\Def(z_j)_{\ell_h(z_j)}$ for $j=0,\dots,N-1$. By
\cref{prop:DPP}, $z_j\in\bm A$ with $\VV(z_j)=\nu-jh$ for $j<N$, so
$\ell_h(z_j)=h$ for $j\le N-2$, $\ell_h(z_{N-1})=\nu-(N-1)h\in(0,h]$,
and $z_N$ is the impact point, in $\Gamma_{\bm A}$. Set
$e_j\Def W(z_j)-\ell_h(z_j)-W(z_{j+1})$, so that $\abs{e_j}\le\delta$
by \cref{eq:residual}. Summing,
$W(z)-\sum_{j<N}\ell_h(z_j)-W(z_N)=\sum_{j<N}e_j$, and
$\sum_{j<N}\ell_h(z_j)=\nu$. Hence
$\abs{W(z)-\nu}\le\abs{W(z_N)}+N\delta\le b+N\delta$. Uniqueness
follows with $\bm A=\mathbb D$, $\delta=b=0$.
\end{proof}

\Cref{thm:certificate} gives a uniform certificate conditional on
verified residual and reached-boundary bounds, with the computable
factor $\ceil{\timehor/h}$. Its proof needs only the finite stopped
chain from each initial point; a residual bound on that chain and its
endpoint also yields a pointwise certificate. Residual samples at
initial points alone do not supply these hypotheses. No derivative or
global bound on $W$ is required. One-sided residual and endpoint
bounds give the corresponding one-sided error bounds.

\begin{proposition}[Exact iteration]\label{prop:iteration}
For every $W_0\in\Bnd$ and $k\in\NN$,
$\Sh^{\,k}[W_0]=\VV$ on $\set{z\in\mathbb D:\VV(z)\le kh}$; in
particular $\Sh^{\,K_h}[W_0]=\VV$ on $\Omega$ for
$K_h\Def\ceil{\timehor/h}$. For an arbitrary bounded $W_0:\Omega\to\RR$
we instead have $\Sh^{\,k+1}[W_0]=\VV$ on
$\set{z\in\mathbb D:\VV(z)\le kh}$ and
$\Sh^{\,K_h+1}[W_0]=\VV$ on $\Omega$, since $\Sh[W_0]\in\Bnd$.
\end{proposition}

\begin{proof}
With the notation of the previous proof, unrolling the iteration from
$z$ with $\nu\le kh$ follows the chain $z_0,\dots,z_N$, $N\le k$,
reaches the impact point $z_N\in\QT$ where every iterate is $0$, and
collects $\sum_{j<N}\ell_h(z_j)=\nu$; later applications of $\Sh$ on
$\QT$ leave the value $0$ unchanged.
\end{proof}

Exact collision information thus propagates backward from the target:
after $k$ sweeps a boundary-compatible initialization is exact on
$\set{\VV\le kh}$; otherwise one additional sweep suffices. This is an identity for
the continuum operator with exact flow and event evaluation; a
numerical implementation with interpolation between nodes incurs the
interpolation error at each sweep.

\paragraph{Extension to the full capped domain.}\label{par:capped-operator}
For numerical states outside $\mathbb D$, define
$\overline\Omega\Def[0,\timehor]\times U$ and the absorbing boundary
$\Gamma_T\Def\QT\cup(\{\timehor\}\times U)$. Set
$\overline S_h[W]=0$ on $\Gamma_T$ and use \cref{E:S-operator}
elsewhere, with the capped value \cref{E:value-function}: advance
until $h$, first contact, or the terminal time, whichever occurs
first. Uniqueness of the prediction flow gives
$\VV(z_\kappa)=\VV(z)-\kappa$ for $0\le\kappa<\VV(z)$, including
grazing and noncollision trajectories. Every chain therefore ends in
$\Gamma_T$ within $K_h$ steps. The fixed-point identity,
\cref{thm:certificate,prop:iteration}, and their proofs extend
verbatim, replacing the reached impact set by the reached part of
$\Gamma_T$ and imposing $W=0$ on $\Gamma_T$ for uniqueness and exact
iteration without an extra sweep. On $\Omega$ this agrees with $\Sh$. This algebraic
extension makes no regularity claim outside $\mathbb D$.

\begin{proposition}[The optional-duration form is not identifiable]
\label{prop:envelope}
Let $\Tdt[W](z)\Def\inf_{0\le\kappa\le\ell_h(z)}\{\kappa+W(z_\kappa)\}$
on $\mathbb D$ and $\Tdt[W]\Def0$ on $\QT$, an optional-duration
alternative to \cref{E:S-operator}. Then:
\begin{enumerate}[label=\tup{\roman*},leftmargin=*,nosep]
\item $\Tdt[W]\le W$ for every $W\in\Bnd$, and $W=\Tdt[W]$ if and only
if $W(z)\le\kappa+W(z_\kappa)$ for all $z\in\mathbb D$,
$0\le\kappa\le\ell_h(z)$; by the semigroup property, if and only if
this holds for all $0\le\kappa\le\VV(z)$. The set $\sub$ of such $W$
does not depend on $h$.
\item $\sub$ is convex, closed under pointwise maxima and minima, and
contains $\lambda\VV$ for every $\lambda\in[0,1]$, in particular $0$.
\item $\VV=\max\sub$: every $W\in\sub$ satisfies $W\le\VV$ on
$\mathbb D$.
\end{enumerate}
\end{proposition}

\begin{proof}
(i) $\kappa=0$ is admissible. For the second equivalence, split
$\kappa\le\VV(z)$ into steps of length at most $h$ along the
trajectory and chain the inequalities using \cref{prop:DPP}. (ii)
The inequality is preserved by convex combinations, maxima and
minima; for $W=\lambda\VV$ it reads $\lambda\nu\le\kappa+\lambda(\nu-\kappa)$,
i.e.\ $(1-\lambda)\kappa\ge0$. (iii) Take $\kappa=\VV(z)$: $W(z)\le\VV(z)+W(z_\nu)
=\VV(z)$.
\end{proof}

Thus the fixed points of the optional-duration operator are exactly
the dynamic-programming subsolutions, a convex family of lower bounds
on the TTC with $\VV$ as its maximum.
Minimizing the residual of $\Tdt$ does not identify the TTC:
$W\equiv0$ attains zero residual and satisfies the boundary condition. The
one-dimensional example $\dot y=-1$, $\mathcal Q=\set{y\le0}$,
$\VV(s,y)=y$ for $0<y<\timehor-s$, shows the contrast: for
$W=\lambda\VV$ and $y>h$, the true error is $\abs{1-\lambda}y$, the
optional-duration residual is $(\lambda-1)_+h$, zero for all
$\lambda\le1$, while the stopped-step residual is $\abs{1-\lambda}h$,
nonzero whenever $W\neq\VV$. The deficiency
is specific to the infimum: removing $\kappa=0$ alone does not cure it
for continuous $W$, since durations arbitrarily close to $0$ remain.

\section{Explicit ReLU Realization}
\label{sec:construction}

In this section, $\|f\|_{L^\infty(E)}$ denotes the pointwise supremum
$\sup_{z\in E}|f(z)|$, including when $E$ has zero Lebesgue measure.
Let $\bm D\subset\mathbb D$ be compact and
$\VV_{\min}(\bm D)\Def\min_{\bm D}\VV>0$. Fix
$h\in(0,\VV_{\min}(\bm D))$, so that $\ell_h=h$ on $\bm D$, and define
the one-step forward tube and the enlarged set
\begin{align*}
\bm D^{h}&\Def\set[\big]{(s+\kappa,\Phi(\kappa;s,\bm y)):
(s,\bm y)\in\bm D,\ \kappa\in[0,h]},\\
\bm D'&\Def\bm D\cup\bm D^{h}.
\end{align*}
$\bm D'$ is compact as a continuous image of $\bm D\times[0,h]$, and
$\bm D'\subset\mathbb D$ by \cref{prop:DPP}; for $z\in\bm D$, the
operator $\Sh$ evaluates its argument at $z_h\in\bm D'$. Since
$\bm D'$ is a compact subset of $[0,\timehor)\times U$ and
$\mathbb D$ is relatively open, there is $r\in(0,\timehor-\max_{\bm D'}s)$,
smaller than the max-norm distance from the $\bm y$-projection of
$\bm D'$ to $\RR^d\setminus U$, such that $\bm K\Def\set{z\in[0,\timehor)\times U:\dist_\infty(z,\bm D')
\le r}$, with $\dist_\infty$ the max-norm distance, is a compact subset
of $\mathbb D$. Let $L\Def L_{\VV}(\bm K)$, let $\bm B\supset\bm K$ be a
closed box with sides $\ell_1,\dots,\ell_m\le\ell_{\bm B}$ whose
$s$-projection lies in $[0,\timehor)$, and set
$\Theta\Def\max\{\timehor,L\ell_{\bm B}\}$.

\begin{theorem}[Optimization-free ReLU realization]\label{thm:relu}
There is a constant $C_m$ depending only on $m$ such that, for every
$\eps\in(0,\Theta)$, there exist grid points
$\set{z_j}_{j\in\mathcal J}\subset\bm K$ and ReLU networks
$\set{\tilde\varphi_j}_{j\in\mathcal J}$, depending only on
$(\bm B,\bm D',r,m,\eps,\timehor,L)$ and not on the values of $\VV$,
such that the network
\begin{equation}\label{eq:relu-net}
\VV_{\theta^*}(z)\Def\sum_{j\in\mathcal J}c_j\,\tilde\varphi_j(z),
\qquad c_j\Def\VV(z_j),
\end{equation}
has depth at most $C_m(\ln(\Theta/\eps)+1)$ and at most
\begin{equation}\label{eq:netsize-new}
\mathcal N(\eps,\bm D')\le C_m\Bigl(\frac{\Theta}{\eps}
+\frac{\ell_{\bm B}}{r}\Bigr)^{m}\Bigl(\ln\frac{\Theta}{\eps}+1\Bigr)
\end{equation}
weights and units in a sparse realization, and satisfies
\begin{align}
\norm{\VV_{\theta^*}-\VV}_{L^\infty(\bm D')}&\le\eps,
\label{eq:Linf-approx}\\
\norm{\VV_{\theta^*}-\Sh[\VV_{\theta^*}]}_{L^\infty(\bm D)}
&\le2\eps .
\label{eq:residual-small}
\end{align}
If the nodal values are replaced by $\tilde c_j$ with
$\abs{\tilde c_j-\VV(z_j)}\le\eta$, then
\cref{eq:Linf-approx,eq:residual-small} hold with $\eps$ replaced by
$\eps+2^{m+1}\eta$.
\end{theorem}

\begin{proof}
\emph{Step 1: interpolation.} Fix $h_0\in(0,r]$ and, on the $i$th
axis of $\bm B$, take the uniform grid of mesh
$h_i\Def\ell_i/\ceil{\ell_i/h_0}\le h_0$, which divides the side
$\ell_i$; the number of nodes on that axis is $\ceil{\ell_i/h_0}+1\le
\ell_i/h_0+2$. Let $\set{z_j}$ be the resulting tensor grid and
$\varphi_j(z)\Def\prod_{i=1}^m\max\{0,1-\abs{z^i-z_j^i}/h_i\}$ the
tensor-product hat functions: $\varphi_j\ge0$, $\sum_j\varphi_j=1$ on
$\bm B$, and $\set{\varphi_j\neq0}\subset\set{\abs{z-z_j}_\infty<h_0}$,
so at most $2^m$ of them are nonzero at any point. Keep in
$\mathcal J$ only the nodes with $\varphi_j\not\equiv0$ on $\bm D'$;
since $h_0\le r$ and the $s$-projection of $\bm B$ lies in
$[0,\timehor)$, these nodes lie in $\bm K$. For $z\in\bm D'$, since $\abs{z-z_j}\le\sqrt m\,h_0$
wherever $\varphi_j(z)\neq0$,
$\abs{\VV(z)-\sum_j\VV(z_j)\varphi_j(z)}\le\sum_j\abs{\VV(z)-\VV(z_j)}
\varphi_j(z)\le L\sqrt m\,h_0$.

\emph{Step 2: ReLU realization.} Each univariate hat is a one-layer
ReLU network. By \cite[Prop.~3 and proof of Thm.~1]{yarotsky2017},
the product of $m$ factors in $[0,1]$ is approximated within
$\varrho\in(0,1]$ by a ReLU network of depth and size
$O(\ln(1/\varrho)+1)$ that vanishes whenever one factor vanishes.
This yields $\tilde\varphi_j$ with $\abs{\tilde\varphi_j-\varphi_j}
\le\varrho$, $\tilde\varphi_j=0$ wherever $\varphi_j=0$, and
$\abs{\tilde\varphi_j}\le2$. As $0\le c_j\le\timehor$ and at most
$2^m$ terms are nonzero,
$\abs{\sum_jc_j(\tilde\varphi_j-\varphi_j)}\le2^m\timehor\varrho$. The
choices $h_0\Def r$ if $L\sqrt m\,r\le\eps/2$ and
$h_0\Def\eps/(2L\sqrt m)$ otherwise, and
$\varrho\Def\min\{1,\eps/(2^{m+1}\timehor)\}$, give
\cref{eq:Linf-approx}. The grid has at most $\prod_i(\ell_i/h_0+2)\le
(\ell_{\bm B}/r+2\sqrt m\,\Theta/\eps+2)^m$ nodes, each costing
$O(\ln(1/\varrho)+1)$ units with $\ln(1/\varrho)\le\ln(\Theta/\eps)
+(m+1)\ln2$; since $\Theta/\eps>1$, this gives \cref{eq:netsize-new}.

\emph{Step 3: residual.} For $z\in\bm D$, \cref{lem:fixed-point} and
\cref{eq:non-expansion}, whose evaluation point $z_h$ lies in
$\bm D'$, give $\abs{\VV_{\theta^*}-\Sh[\VV_{\theta^*}]}(z)\le
\abs{\VV_{\theta^*}-\VV}(z)+\abs{\VV_{\theta^*}-\VV}(z_h)\le2\eps$.

\emph{Step 4: perturbed nodes.} The additional error is
$\abs{\sum_j(\tilde c_j-c_j)\tilde\varphi_j}\le2^m\cdot2\eta$, and
Step~3 applies unchanged.
\end{proof}

\begin{remark}[Scope of \cref{thm:relu}]\label{rem:scope}
After the geometric and regularity parameters $(\bm B,\bm D',r,L)$
are fixed, the only coefficients carrying value data are the nodal
TTC values $c_j$, obtained offline by direct simulation or by the
sweeps of \cref{prop:iteration}; their numerical errors enter through
$\eta$. Once valid geometric and regularity bounds are supplied,
the theorem gives a priori accuracy, size and depth estimates, and
no loss is minimized. The bound
\cref{eq:netsize-new} is an a priori complexity estimate with an
unspecified $C_m$, exponential in $m$ as expected for Lipschitz
functions \cite{yarotsky2017}. The approximation is local to a compact subset of the
transversal collision domain, separated from grazing, from the
collision set and from the terminal horizon; uniform bounds need
not persist as these separations are lost. The local residual bound
\cref{eq:residual-small} alone does not supply the full-chain
hypotheses of \cref{thm:certificate}: a compact interior set
$\bm D$ cannot contain complete pre-contact trajectories.
Using that certificate additionally requires residual bounds along
the remaining chain and bounds at the reached boundary points.
\end{remark}

\section{Numerical Illustration}
\label{sec:numerics}
\textit{The codes and detailed results can be found in \cite{companion}}.
We use \cref{E:new-accel} with $\alpha=2$, $\timehor=1$ and
$\sa_{\mathsf L}=0.5,-1.2,0.8,-0.4,0$ on the fifths of $[0,1]$.
Initial $(s,\Xi,v_{\mathsf F},\zeta_{\mathsf L})$ are independent
uniform draws on $[0,1)\times(0,1]\times[0,v_{\max}]\times[0.25,1.75]$,
where $v_{\max}=2\tanh2$ and $\Lambda=\zeta_{\mathsf L}-v_{\mathsf F}$.
This is a sampling cap; the global supremum of $\mathcal V$ is
$1+\tanh2$. For these initial conditions the leader speed stays in
$[0.01,1.91]$, the pre-contact gap is below $4$, and hence
$0\le v_{\mathsf F}\le v_{\max}$ is preserved over this horizon.
Reference TTCs use RK4 with maximum step $10^{-3}$, split at forcing
jumps, stopping times and the horizon, with $40$ event bisections.
Initial contact and terminal states are absorbed immediately; a
closing-to-opening step also checks its interior minimum. Independent
piecewise DOP853 comparisons and step refinements on $120$ selected
cases give TTC differences below $10^{-10}$. These finite checks do
not give a uniform error enclosure or certify grazing classifications.

\textbf{Learned candidates and independent evaluation.}
The two ReLU networks have three width-$400$ hidden layers ($1{,}200$
ReLU units, $323{,}201$ parameters). Both use $10^6$ training and
$10^5$ validation states, with initial times rounded to $10^{-3}$;
$44{,}959$ training impact states supply zero boundary targets.
$\VV_\theta$ minimizes MSE, while $\VV_\theta^S$ also penalizes the
mean squared stopped-step residual with $h=0.05$ and weight $1$.
Training uses Adam, $100$ epochs, batch size $4096$, and initial
learning rate $10^{-3}$ with cosine decay. The trained checkpoints
are held fixed for the following evaluation.

\begin{figure*}[t]
\centering
\includegraphics[width=\textwidth]{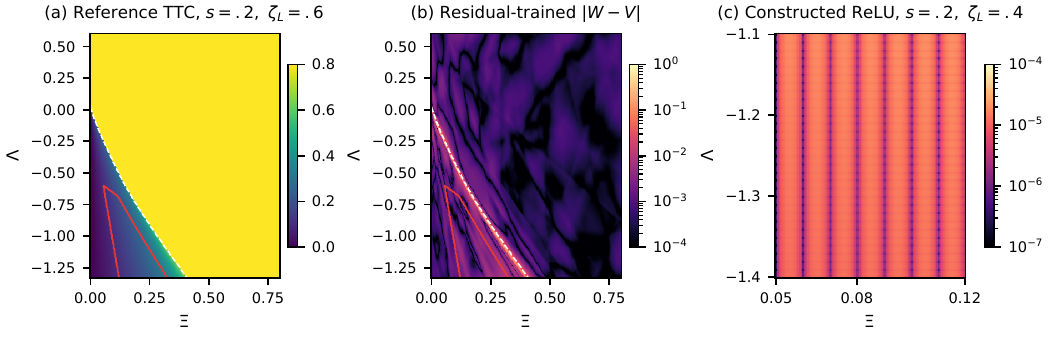}
\caption{(a),(b) Physical slice $s=0.2$, $\zeta_{\mathsf L}=0.6$:
reference TTC and error of the residual-trained candidate. White
dash: sampled contact/noncontact interface; red: the selection used
for $\widehat D$. (c) Actual constructed ReLU error on a different,
fully covered slice of $B_c$, $s=0.2$, $\zeta_{\mathsf L}=0.4$,
$h_0=0.01$. The axes and error color scale in (c) differ from (b).}
\label{fig:errors}
\end{figure*}

A fresh test of $300{,}000$ states uses continuous, unrounded times
and seed $20261103$; $4.51\%$ collide. Its selected subset
$\widehat D$ has $2{,}725$ points with $\VV\ge0.1$,
$\Lambda(\VV;z)\le-0.5$, and $\VV\le1-s-0.2$.
For $(\VV_\theta,\VV_\theta^S)$, maximum errors on $\widehat D$ are
$(0.00924,0.01353)$ and RMSEs are $(0.002289,0.002206)$.
Residual training modestly improves RMSE but worsens the maximum
error. Full-sample maxima are $(0.8881,0.9333)$, with large errors
near the contact/noncontact interface; the regularity and
approximation results make no uniform claim there
(\cref{fig:errors}(a),(b)).

\begin{figure}[t]
\centering
\includegraphics[width=\linewidth]{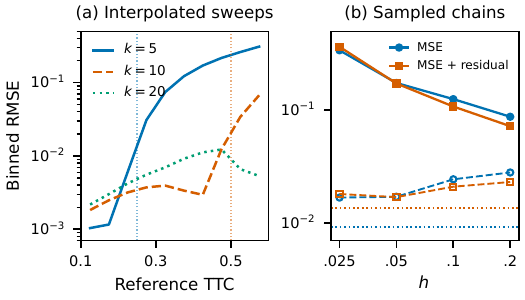}
\caption{(a) Binned RMSE for interpolated stopped-step sweeps from
$W_0=0$ on the fine grid's selected nodes; vertical lines mark
$\VV=kh$ for $k=5,10$. (b) Finite-sample envelopes (solid), chain
residuals (dashed), and observed maximum errors (dotted), over
$\widehat D$. These curves do not establish uniform certification.}
\label{fig:iter}
\end{figure}

\textbf{Sampled stopped chains.}
For each initial point of $\widehat D$, we evaluate every node of
its complete stopped chain and its reached impact state, for
$h\in\{0.025,0.05,0.1,0.2\}$. Let $\widehat\delta_h$ and
$\widehat b$ be the maxima of the absolute residual and boundary
value over these finite collections. At $h=0.05$ they are
$(0.01712,0.016996)$ and $(0.003318,0.002870)$, respectively.
Every sampled error lies below
$\ceil{\VV/h}\widehat\delta_h+\widehat b$, whose maxima are
$(0.17452,0.17283)$. The largest tested algebraic telescoping defect
is below $8\times10^{-8}$. These are empirical finite-chain checks
of \cref{thm:certificate}, not certificates on a continuous domain;
flow, event and roundoff errors are not uniformly bounded.
The envelopes decrease with $h$ in this experiment
(\cref{fig:iter}(b)). For comparison, $W=\lambda\VV$ with
$\lambda=0.5$ has zero optional-duration residual analytically on any exact
chain; its maximum error over the sampled initial points is $0.2432$; its stopped-step
residual is $0.025$ at $h=0.05$, as \cref{prop:envelope} predicts.

\textbf{Complete four-dimensional construction.}
Independently of all test draws, fix the boxes, in coordinates
$(s,\Xi,\Lambda,\zeta_{\mathsf L})$,
\begin{align*}
B_c&=[.15,.27]\times[.05,.12]\times[-1.4,-1.1]\times[.27,.52],\\
D_c&=[.15,.25]\times[.08,.12]\times[-1.4,-1.2]\times[.30,.50].
\end{align*}
The whole box $B_c$ lies in $\mathbb D$. Indeed, until contact,
$\zeta_{\mathsf L}\in[.03,.76]$ and the closing speed
$w=-\Lambda$ satisfies $w\le1.4$ and
$w\ge2.26e^{-2\kappa}-1.16$. At
$T_*\Def\tfrac12\log(2.26/1.46)\approx0.218464$ this lower bound is
$c=.3$, while its integral is $0.4-1.16T_*>.12$.
Thus $.05/1.4\le\VV<T_*$ and impact is uniformly transversal.
The complete $h_c=.02$ tube from $D_c$ lies in $B_c$:
$\Xi\ge.052$, $1.1<w\le1.4$, and
$\zeta_{\mathsf L}\in[.276,.516]$ along that tube.

An analytical Lipschitz constant on $B_c$ follows by subtracting
the integrated follower equations. Write
$\ell=\operatorname{sech}^2(1.88)$ and
$b=(1-e^{-2T_*})/2$. For two initial states, before their earlier
contact, the gap difference is bounded by
\[
\frac{|\Delta\Xi_0|+b|\Delta\Lambda_0|+(T_*-b)|\Delta\zeta_{\mathsf L,0}|
       +2(1.2)T_*|\Delta s|}{1-\ell T_*^2}.
\]
Here the shifted leader integrals differ by at most
$2\|\sa_{\mathsf L}\|_\infty|\Delta s|$; no derivative of the input
is used. Dividing by the closing-speed bound $c$ gives
\[
L=\frac{\sqrt{(2.4T_*)^2+1+b^2+(T_*-b)^2}}
           {.3(1-\ell T_*^2)}<3.82845.
\]
We retain every node of a tensor grid on $B_c$, giving a complete-box
version of \cref{sec:construction}; $h_0=.02,.01,.005$ use
$7{,}840$, $83{,}824$, and $1{,}166{,}625$ nodes. No query determines
which nodes are retained. The actual ReLU products are evaluated:
if $q_k$ is the dyadic linear interpolant of $x^2$, the multiplier
$[q_k((x+y)/2)-q_k(|x-y|/2)]_{[0,1]}$ is realized with ReLUs and has
error at most $e_k=4^{-(k+1)}$, preserving zero factors exactly.
For $k=6,7,8$, the exact-node uniform error bound is
\begin{equation}\label{eq:numeric-construction-bound}
E_h=\frac L2\Bigl(\sum_{i=1}^4h_i^2\Bigr)^{1/2}+48T_*e_k.
\end{equation}
The first term follows from the variance of the multilinear
barycentric weights; the second uses at most $16$ active hats and
three multipliers. Nodal errors bounded by $\eta$ add at most
$(1+48e_k)\eta$, and the local residual on $D_c$ is at most twice
the resulting bound. The numerical nodal error is not uniformly
certified; $116$ independent reference comparisons and $100$ nodal
checks per grid give discrepancies below $10^{-10}$.

\begin{table}[t]
\centering\scriptsize
\caption{Common independent $10{,}000$-point test on $D_c$.
$N_\sigma$: ReLU gates in an unshared feedforward graph with affine
skip connections, not parameters. $E_h$: analytical bound with exact
nodal values; computed nodal errors require the additional term
in the text. Learned models are the unchanged broad-domain models.}
\label{tab:comparison}
\setlength{\tabcolsep}{3.2pt}
\begin{tabular}{lrrrr}
\hline
Candidate & $N_\sigma$ & max. error & RMSE & $E_h$\\
\hline
$\VV_\theta$ & $1{,}200$ & $2.75\,10^{-3}$ & $9.13\,10^{-4}$ & --\\
$\VV_\theta^S$ & $1{,}200$ & $3.93\,10^{-3}$ & $1.24\,10^{-3}$ & --\\
$h_0=.02$ & $1.03\,10^6$ & $8.57\,10^{-4}$ & $1.68\,10^{-4}$ & $.07418$\\
$h_0=.01$ & $1.26\,10^7$ & $2.11\,10^{-5}$ & $9.29\,10^{-6}$ & $.03844$\\
$h_0=.005$ & $1.96\,10^8$ & $5.43\,10^{-6}$ & $2.31\,10^{-6}$ & $.01918$\\
\hline
\end{tabular}
\end{table}

\Cref{tab:comparison} compares all candidates on the same $10{,}000$
independent $D_c$ states. A separate $10{,}000$-point test throughout
$B_c$ gives constructed maxima
$9.82\,10^{-4}$, $2.89\,10^{-5}$, and $9.49\,10^{-6}$;
the corresponding sampled local residuals on $D_c$ are
$8.07\,10^{-4}$, $3.03\,10^{-5}$, and $7.50\,10^{-6}$.
The coarse time grid misses the forcing jump at $.2$, whereas the
two finer grids contain it. The construction is accurate here but
large; its cost is not comparable to training over the larger
sampling domain. Evaluation uses compact support to sum the same
ReLU formula, and is not a dense-network latency measurement.

\textbf{Interpolated sweeps.}
We use the capped-domain extension $\overline S_h$ in
\cref{sec:stopped}, with contact and the terminal plane absorbing.
Let $A(s)=\int_0^s\sa_{\mathsf L}(t)\,dt$, and grid the coordinates
$(s,q,v_{\mathsf F},\omega)$ with
$q=\Xi/(1+2s)$, $\omega=\zeta_{\mathsf L}-A(s)$ on
$[0,1]^2\times[0,1+\tanh2]\times[.15,1.89]$.
This enlarged computational box contains the sampling domain and is
invariant until stopping: $\omega$ is constant, the follower-speed
interval is invariant, and
$\dot q=(\Lambda-2q)/(1+2s)<0$ at $q=1$ since
$\zeta_{\mathsf L}\le1.99$. No escaped state is projected back.
With $h=.05$, grid shapes $21\times31\times17\times17$ and
$21\times61\times33\times33$ have $188{,}139$ and $1{,}395{,}009$
nodes. Only $(q,v_{\mathsf F})$ need bilinear interpolation at the
regular landing points. On the same $288$ independently drawn
selected states, the final maximum error decreases from $.04873$ to
$.02584$ and RMSE from $.01735$ to $.00601$ under refinement.
Initializations $W_0=0$ and $W_0=.5$ reach the same discrete values
by sweeps $20$ and $21$, respectively. Spatial interpolation error
remains after finite propagation; \cref{prop:iteration} is a result
for the exact continuum operator. These experiments evaluate a
prescribed prediction map, and do not test closed-loop safety.

\section{Conclusion}
\label{sec:conclusion}

Dynamic TTC under a prescribed Carath\'eodory prediction flow is
locally Lipschitz on the transversal collision domain and is uniquely
characterized by a derivative-free stopped-step Bellman equation.
Complete finite chains give a two-sided residual error bound and
finite recovery by the exact operator; optional-duration infimum
residuals alone do not identify the value function. A local ReLU
construction gives explicit approximation and nodal-perturbation
bounds. The computations illustrate these mechanisms with independent
tests, a fully covered four-dimensional construction, and grid
refinement. Uniform numerical certification still requires validated
flow, event, nodal and roundoff error bounds. Transferring TTC
approximation guarantees to the closed loop remains a separate
control problem.

\appendix
\section*{Uniform boundedness up to contact}

\begin{lemma}\label{lem:uniform-bound}
Let \ref{A1}--\ref{A2} hold and let $\bm C_0\subset P$ be compact.
There is a compact $\bm C_1\subset P$ such that every maximal solution
in $P$ from a datum in $[0,\timehor]\times\bm C_0$ takes its values in
$\bm C_1$; in particular $\Phi(\kappa;s,\bm y)\in\bm C_1$ for all
$(s,\bm y)\in[0,\timehor]\times\bm C_0$, $0\le\kappa\le\varsigma(s,\bm y)$.
\end{lemma}

\begin{proof}
Two consequences of \ref{A1}. \emph{(F1)} If $x,x'$ are solutions in
$P$ on $[a,b]$ with values in a compact $\bm K\subset P$, then
$\abs{x(t)-x'(t)}\le e^{L_{\bm K}(t-a)}\abs{x(a)-x'(a)}$ on $[a,b]$
and $\abs{x(t)-x(t')}\le M_{\bm K}\abs{t-t'}$. \emph{(F2)} If $x_n$
are solutions in $P$ on $[b_n,b_n+\theta]$ with $b_n\to b$, values in
a compact $\bm K\subset P$, and $x_n(b_n)\to\bar{\bm y}$, then a
subsequence converges uniformly on $[b,b+\theta/2]$, after extension
by the constant $x_n(b_n)$ on $[b,b_n]$ when $b_n>b$, to a solution in
$P$ on $[b,b+\theta/2]$ from $(b,\bar{\bm y})$: the extended functions
are $M_{\bm K}$-Lipschitz and bounded, so Arzel\`a--Ascoli gives a
uniformly convergent subsequence with limit $\hat x$, with values in
the closed set $\bm K$; for $t\in(b,b+\theta/2]$ and $n$ large,
$x_n(t)=x_n(b_n)+\int_{b_n}^t\mathsf G(r,x_n(r))\,dr$, and since
$\mathsf G(r,\cdot)$ is Lipschitz on $\bm K$ for a.e.\ $r$ and bounded
by $M_{\bm K}$, dominated convergence gives
$\hat x(t)=\bar{\bm y}+\int_b^t\mathsf G(r,\hat x(r))\,dr$.

Suppose the lemma fails. Let $P_1\subset P_2\subset\cdots$ be compact
subsets of $P$ with $\bigcup_kP_k=P$ and $P_k$ contained in the
interior of $P_{k+1}$ relative to $P$ \tup{such an exhaustion exists
since $P$ is closed in the open set $U$}. Then for every $k$ some
maximal solution $x_k$ in $P$ from $(s_k,\bm y_k)\in[0,\timehor]\times
\bm C_0$, on $[s_k,s_k+\sigma_k]$, satisfies $x_k(t_k)\notin P_k$ for
some $t_k$. Passing to a subsequence, $(s_k,\bm y_k)\to(s^*,\bm y^*)$.
Let $x^*$ be the maximal solution in $P$ from $(s^*,\bm y^*)$, on
$[s^*,e^*]$; by \ref{A2} its graph is compact. Fix
$\eta\in(0,\dist(x^*([s^*,e^*]),\partial U))$, with
$\dist(\cdot,\varnothing)=\infty$, and set
$\bm K\Def\set{\bm y\in P:\dist(\bm y,x^*([s^*,e^*]))\le\eta}$, a
compact subset of $P$ since the closed $\eta$-neighborhood of the
graph is a compact subset of $U$ and $P$ is closed relative to $U$;
let $M\Def\max\{1,M_{\bm K}\}$, $L\Def L_{\bm K}$. We show that $x_k$ takes
values in $\bm K$ for all large $k$.

Let $\delta_k\Def\abs{\bm y_k-\bm y^*}+M\abs{s_k-s^*}\to0$ and
$T_k\Def\sup\set{t\in[s_k,s_k+\sigma_k]:x_k([s_k,t])\subset\bm K}$;
for $k$ large $\bm y_k\in\bm K$, and $x_k([s_k,T_k])\subset\bm K$ by
continuity. Assume $T_k<s_k+\sigma_k$ for infinitely many $k$; we show
that $x_k(T_k)$ then lies within $\eta/2$ of $x^*([s^*,e^*])$, so that
$x_k$, being $M$-Lipschitz while in $\bm K$, remains in $\bm K$ on
$[T_k,\min\{T_k+\eta/(2M),s_k+\sigma_k\}]$, contradicting the
definition of $T_k$.

\emph{Step 1: the overlap.} Let $a_k\Def\max\{s_k,s^*\}$ and suppose
$a_k\le\min\{T_k,e^*\}$. On $[\min\{s_k,s^*\},a_k]$ only one solution
moves, inside $\bm K$, so by (F1) $\abs{x_k(a_k)-x^*(a_k)}\le\delta_k$,
and on $[a_k,\min\{T_k,e^*\}]$ both lie in $\bm K$, so
\begin{equation}\label{eq:overlap}
\abs{x_k(t)-x^*(t)}\le e^{L\timehor}\delta_k,\qquad
t\in[a_k,\min\{T_k,e^*\}].
\end{equation}
If $T_k\le e^*$ and $a_k\le T_k$, \cref{eq:overlap} at $t=T_k$ gives
the announced contradiction for $k$ large. If $T_k<a_k$ (only if
$s_k<s^*$), then $T_k-s_k<s^*-s_k$ and $\abs{x_k(T_k)-\bm y^*}\le
\delta_k$, again a contradiction. Hence $T_k>e^*$ for $k$ large.

\emph{Step 2: beyond $e^*$.} Since $T_k\le\timehor$, $e^*<\timehor$,
so $x^*$ ends at $e^*$ and admits no continuation in $P$. Let
$b_k\Def\max\{s_k,e^*\}$ and $\theta_k\Def T_k-b_k>0$. By
\cref{eq:overlap} at $t=e^*$ when $s_k\le e^*$, and by
$x_k(b_k)=\bm y_k$ otherwise (then $e^*=s^*$), $x_k(b_k)\to x^*(e^*)$.
If $\theta_k\ge\theta>0$ along a subsequence, (F2) applied to the
restrictions of $x_k$ to $[b_k,b_k+\theta]$ yields a solution in $P$
on $[e^*,e^*+\theta/2]$ from $(e^*,x^*(e^*))$, whose concatenation
with $x^*$ contradicts maximality. Hence $\theta_k\to0$ and
$\abs{x_k(T_k)-x^*(e^*)}\le\abs{x_k(b_k)-x^*(e^*)}+M\theta_k\le\eta/2$
for $k$ large, the announced contradiction.

Therefore $x_k$ takes values in $\bm K$ for all large $k$. Since
$\bm K$ is a compact subset of $P$, it is contained in $P_k$ for $k$
large, contradicting $x_k(t_k)\notin P_k$. Hence there is $k_0$ such
that every maximal solution from $[0,\timehor]\times\bm C_0$ takes
values in $P_{k_0}$, and $\bm C_1\Def P_{k_0}$ is the required set.
\end{proof}

\bibliographystyle{IEEEtran}
\bibliography{reference_revised}

@inproceedings{matin2020nonlinear,
  author = {{Nick Zinat Matin}, Hossein and Sowers, Richard B.},
  title = {Nonlinear optimal velocity car following dynamics ({I}): Approximation in presence of deterministic and stochastic perturbations},
  booktitle = {American Control Conference},
  pages = {410--415}, year = {2020}
}

@inproceedings{matin2020nonlinear2,
  author = {{Nick Zinat Matin}, Hossein and Sowers, Richard B.},
  title = {Nonlinear optimal velocity car following dynamics ({II}): Rate of convergence in the presence of fast perturbation},
  booktitle = {American Control Conference},
  pages = {416--421}, year = {2020}
}

@article{matin2024analytical,
  author = {{Nick Zinat Matin}, Hossein and Yeo, Yuneil and Gong, Xiaoqian and {Delle Monache}, Maria Laura},
  title = {On the analytical properties of a nonlinear microscopic dynamical model for connected and automated vehicles},
  journal = {IEEE Control Systems Letters},
  volume = {8}, pages = {1607--1612}, year = {2024},
  doi = {10.1109/LCSYS.2024.3412156}
}

@misc{matin2026TTCandHJB,
  author = {{Nick Zinat Matin}, Hossein and Yeo, Yuneil and {Delle Monache}, Maria Laura},
  title = {Regularity of the time-to-collision map and convergence of learned {Hamilton--Jacobi} solutions},
  howpublished = {Preprint}, year = {2026},
  note = {\href{https://github.com/hosein-matin/Papers/blob/edb3eb0ea9c3333bbcf02887157e53c0b3ce3a6f/Regularity\%20of\%20Time-to-Collision\%20Map\%20and\%20Convergence\%20of\%20Learned\%20HJ\%20Solution\%20.pdf}{Archived snapshot, commit \texttt{edb3eb0}, Oct.\ 27, 2026}}
}

@techreport{hayward1971near,
  author = {Hayward, J.},
  title = {Near misses as a measure of safety at urban intersections},
  institution = {Pennsylvania Transportation and Traffic Safety Center},
  year = {1971}
}

@InProceedings{matin2025second,
author = {Matin, Hossein Nick Zinat and Yeo, Yuneil and Ngo, Amelie Ju-Kang and Paiva, Antonio R and Utke, Jean and Monache, Maria Laura Delle},
title = "Second-Order Time to Collision With Non-Static Acceleration",
booktitle = "Traffic and Autonomy",
year = "2026",
publisher = "Springer Nature Switzerland",
address = "Cham",
pages = "37--83",
isbn = "978-3-032-24629-5"
}

@article{lee1976theory,
  author = {Lee, D. N.},
  title = {A theory of visual control of braking based on information about time-to-collision},
  journal = {Perception},
  volume = {5}, number = {4}, pages = {437--459}, year = {1976}
}

@article{mitchell2005time,
  author = {Mitchell, Ian M. and Bayen, Alexandre M. and Tomlin, Claire J.},
  title = {A time-dependent {Hamilton--Jacobi} formulation of reachable sets for continuous dynamic games},
  journal = {IEEE Transactions on Automatic Control},
  volume = {50}, number = {7}, pages = {947--957}, year = {2005}
}

@inproceedings{bansal2017hamilton,
  author = {Bansal, Somil and Chen, Mo and Herbert, Sylvia and Tomlin, Claire J.},
  title = {{Hamilton--Jacobi} reachability: A brief overview and recent advances},
  booktitle = {IEEE Conference on Decision and Control},
  pages = {2242--2253}, year = {2017}
}

@article{raissi2019,
  author = {Raissi, Maziar and Perdikaris, Paris and Karniadakis, George Em},
  title = {Physics-informed neural networks: A deep learning framework for solving forward and inverse problems involving nonlinear partial differential equations},
  journal = {Journal of Computational Physics},
  volume = {378}, pages = {686--707}, year = {2019}
}

@article{veliov1997lipschitz,
  author = {Veliov, V. M.},
  title = {Lipschitz continuity of the value function in optimal control},
  journal = {Journal of Optimization Theory and Applications},
  volume = {94}, pages = {335--363}, year = {1997},
  doi = {10.1023/A:1022683628650}
}

@article{motta2014asymptotic,
  author = {Motta, Monica and Sartori, Caterina},
  title = {On asymptotic exit-time control problems lacking coercivity},
  journal = {ESAIM: Control, Optimisation and Calculus of Variations},
  volume = {20}, number = {4}, pages = {957--982}, year = {2014}
}

@article{manne1960linear,
  author = {Manne, Alan S.},
  title = {Linear programming and sequential decisions},
  journal = {Management Science},
  volume = {6}, number = {3}, pages = {259--267}, year = {1960}
}

@article{defarias2003linear,
  author = {{de Farias}, Daniela P. and {Van Roy}, Benjamin},
  title = {The linear programming approach to approximate dynamic programming},
  journal = {Operations Research},
  volume = {51}, number = {6}, pages = {850--865}, year = {2003}
}

@article{bardi1990approximation,
  author = {Bardi, Martino and Falcone, Maurizio},
  title = {An approximation scheme for the minimum time function},
  journal = {SIAM Journal on Control and Optimization},
  volume = {28}, number = {4}, pages = {950--965}, year = {1990}
}

@article{yarotsky2017,
  author = {Yarotsky, Dmitry},
  title = {Error bounds for approximations with deep {ReLU} networks},
  journal = {Neural Networks},
  volume = {94}, pages = {103--114}, year = {2017}
}

@article{he2020relu,
  author = {He, Juncai and Li, Lin and Xu, Jinchao and Zheng, Chunyue},
  title = {{ReLU} deep neural networks and linear finite elements},
  journal = {Journal of Computational Mathematics},
  volume = {38}, number = {3}, pages = {502--527}, year = {2020},
  doi = {10.4208/jcm.1901-m2018-0160}
}

@article{oustry2023minimal,
  author = {Oustry, Antoine and Tacchi, Matteo},
  title = {Minimal time nonlinear control via semi-infinite programming},
  journal = {SIAM Journal on Control and Optimization},
  volume = {64}, number = {2}, pages = {698--726}, year = {2026},
  doi = {10.1137/23M1582552}
}

@article{liu2026verifiable,
  author = {Liu, Jun},
  title = {Verifiable error bounds for physics-informed neural network solutions of {Lyapunov} and {Hamilton--Jacobi--Bellman} equations},
  journal = {IEEE Control Systems Letters},
  volume = {10}, pages = {403--408}, year = {2026},
  doi = {10.1109/LCSYS.2026.3696322}
}

@article{nick2022near,
  author = {{Nick Zinat Matin}, Hossein and Sowers, Richard B.},
  title = {Near-collision dynamics in a noisy car-following model},
  journal = {SIAM Journal on Applied Mathematics},
  volume = {82}, number = {6}, pages = {2080--2110}, year = {2022}
}

@article{bando1995dynamical,
  author = {Bando, Masako and Hasebe, Katsuya and Nakayama, Akihiro and Shibata, Akihiro and Sugiyama, Yuki},
  title = {Dynamical model of traffic congestion and numerical simulation},
  journal = {Physical Review E},
  volume = {51}, number = {2}, pages = {1035--1042}, year = {1995}
}

@book{filippov1988differential,
  author = {Filippov, A. F.},
  title = {Differential Equations with Discontinuous Righthand Sides},
  publisher = {Kluwer Academic Publishers}, address = {Dordrecht}, year = {1988}
}

@misc{companion,
  author = {H. {Nick Zinat Matin}},
  title = {Derivative free Bellman qquation for TTC with two-sided certificate: code and numerical companion},
  howpublished = {\url{https://github.com/hosein-matin/TTC_HJB_Weak}},
  year = {2026}
}

\end{document}